\documentclass{amsart}
\usepackage{amscd}
\usepackage{amsmath}
\usepackage{amsxtra}
\usepackage{amsfonts}
\usepackage{amssymb}

\newtheorem{theorem}{Theorem}[section]
\newtheorem{corollary}[theorem]{Corollary}
\newtheorem{lemma}[theorem]{Lemma}
\newtheorem{proposition}[theorem]{Proposition}
\theoremstyle{definition}
\newtheorem{definition}[theorem]{Definition}
\newtheorem{remark}[theorem]{Remark}

\theoremstyle{remark}

\renewcommand{\theclaim}{\textup{\theclaim}}

\newtheorem*{acknowledgements}{Acknowledgements}

\numberwithin{equation}{section}

\def\openone%{\hbox{\upshape \small1\kern-3.3pt\normalsize1}}

{\mathchoice

{\hbox{\upshape \small1\kern-3.3pt\normalsize1}}

{\hbox{\upshape \small1\kern-3.3pt\normalsize1}}

{\hbox{\upshape \tiny1\kern-2.3pt\SMALL1}}

{\hbox{\upshape \Tiny1\kern-2pt\tiny1}}}

\makeatletter

\newbox\ipbox

\newcommand{\ip}[2]{\left\langle #1\,|\,#2\right\rangle}

\newcommand{\diracb}[1]{\left\langle #1\mathrel{\mathchoice

{\setbox\ipbox=\hbox{$\displaystyle \left\langle\mathstrut 
#1\right.$}

\vrule height\ht\ipbox width0.25pt depth\dp\ipbox}

{\setbox\ipbox=\hbox{$\textstyle \left\langle\mathstrut 
#1\right.$}

\vrule height\ht\ipbox width0.25pt depth\dp\ipbox}

{\setbox\ipbox=\hbox{$\scriptstyle \left\langle\mathstrut 
#1\right.$}

\vrule height\ht\ipbox width0.25pt depth\dp\ipbox}

{\setbox\ipbox=\hbox{$\scriptscriptstyle \left\langle\mathstrut 
#1\right.$}

\vrule height\ht\ipbox width0.25pt depth\dp\ipbox}

}\right. }

\newcommand{\dirack}[1]{\left. \mathrel{\mathchoice

{\setbox\ipbox=\hbox{$\displaystyle \left.\mathstrut 
#1\right\rangle$}

\vrule height\ht\ipbox width0.25pt depth\dp\ipbox}

{\setbox\ipbox=\hbox{$\textstyle \left.\mathstrut 
#1\right\rangle$}

\vrule height\ht\ipbox width0.25pt depth\dp\ipbox}

{\setbox\ipbox=\hbox{$\scriptstyle \left.\mathstrut 
#1\right\rangle$}

\vrule height\ht\ipbox width0.25pt depth\dp\ipbox}

{\setbox\ipbox=\hbox{$\scriptscriptstyle \left.\mathstrut 
#1\right\rangle$}

\vrule height\ht\ipbox width0.25pt depth\dp\ipbox}

} #1\right\rangle}

\newcommand{\bz}{\mathbb{Z}}
\newcommand{\br}{\mathbb{R}}
\newcommand{\bc}{\mathbb{C}}

\newcommand{\bn}{\mathbb{N}}

\newcommand{\xir}{X_\infty(r)}
\usepackage{graphicx}

\usepackage{diagrams}

\begin{document}
\title[Measures on projective limit spaces]{Operators, martingales, and measures on projective limit spaces}
\author{Dorin Ervin Dutkay and Palle E.T. Jorgensen}
\address{Department of Mathematics\\
The University of Iowa\\
14 MacLean Hall\\
Iowa City, IA 52242-1419\\
U.S.A.\\} \email{Dorin Ervin Dutkay: ddutkay@math.uiowa.edu}\ 
\email{Palle E.T. Jorgensen: jorgen@math.uiowa.edu}
\subjclass[2000]{42C40, 42A16, 42A65, 43A65, 46G15, 47D07, 60G18}
\keywords{Measures, projective limits, transfer operator, 
martingale, fixed point} 
\dedicatory{Dedicated to the memory of Shizuo Kakutani.}
\begin{abstract}
Let $X$ be a compact Hausdorff space. We study finite-to-one 
mappings $r\colon  X\rightarrow X$, onto $X$, and measures on the 
corresponding projective limit space $X_\infty(r)$. We show that 
the invariant measures on $X_\infty(r)$ correspond in a 
one-to-one fashion to measures on $X$ which satisfy two 
identities. Moreover, we identify those special measures on 
$X_\infty(r)$ which are associated via our correspondence with a 
function $V$ on $X$, a Ruelle transfer operator $R_V$, and an 
equilibrium measure $\mu_V$ on $X$.
\end{abstract}
\maketitle 

\tableofcontents
\section{\label{intr}Introduction}

This paper is motivated by our desire to apply wavelet methods to some
nonlinear problems in symbolic and complex dynamics. Recent research by many
authors (see, e.g., \cite{AST04,ALTW04,BJMP04}) on iterated function systems 
with
affine scaling has suggested that the scope of the multiresolution method
is wider than the more traditional wavelet context. 

Let $r$ be an endomorphism in a compact metric space $X$ (for example the Julia
set of a given rational map $w = r(z)$), and suppose $r$ is onto $X$ and
finite-to-one. Form a projective space $P=P(X,r)$ such that $r$ induces an
automorphism $a=a(r)$ of $P(X,r)$. Let $V$ be a Borel function on $X$ (naturally
extended to a function on $P$). Generalizing the more traditional approach to
scaling functions, we give in theorem \ref{th5_1} a complete classification of
measures on $P(X,r)$ which are quasi-invariant under $a(r)$ and have
Radon-Nikodym derivative equal to $V$. Our analysis of the quasi-invariant
measures is based on certain Hilbert spaces of martingales (theorem \ref{th3_3}) and
on a transfer operator (section \ref{projspec}) studied first by David Ruelle. In
theorem \ref{th7_2} we give a characterization of the extreme points in the set of
$V$-quasi-invariant probability measures.
\par
A basic tool in stochastic processes (from probability theory) 
involves a construction on a ``large'' projective space $X_\infty$, 
based on some marginal measure on some coordinate space $X$. In 
this paper, we consider a special case of this: Our projective 
limit space $X_\infty$ will be constructed from a finite 
branching process.
\par
Our starting point is a finite-to-one mapping $r\colon X\rightarrow X$. 
We will assume that $r$ is onto, but not invertible. Hence, each 
point $x$ in $X$ is finitely covered, and we get a corresponding 
random walk down the iterated branches of powers of $r$, i.e., 
$$r^n:=\underbrace{r\circ r\circ\dots \circ r}_{n\mbox 
{ times}},\,\mbox{ composition of }r\mbox{ with itself }n\mbox{ 
times.}$$ We prescribe probabilities via a fixed function 
$V\colon X\rightarrow[0,\infty)$.
\par
In \cite{DuJo04}, we showed that this setup arises naturally in 
connection with the analysis of wavelets, measures on Julia sets 
for iterations of rational mappings,
$$r(z)=\frac{p_1(z)}{p_2(z)}\,
\mbox{where }p_1\mbox{ and }p_2\mbox{ are polynomials without 
common factors},$$$$\mbox{ and }X\mbox{ is the corresponding 
Julia set;}$$ and for subshift dynamical systems.
\par
Let $A$ be a $k\times k$ matrix with entries in $\{0,1\}$. 
Suppose every column in $A$ contains an entry $1$.
\par
Set
$$X(A):=\left\{(\xi_i)_{i\in\bn}\in\prod_{\bn}\{1,\dots ,k\}\Bigm|A(\xi_i,\xi_{i+1})=1\right\}$$
and
$$r_A(\xi_1,\xi_2,\dots )=(\xi_2,\xi_3,\dots )\mbox{ for }\xi\in X(A).$$
Then $r_A$ is a subshift, and the pair $(X(A),r_A)$ satisfies our 
conditions.
\par
It is known \cite{Rue89} that, for each $A$, as described, the 
corresponding system $r_A\colon X(A)\rightarrow X(A)$ has a unique 
strongly $r_A$-invariant probability measure, $\rho_A$, i.e., a 
probability measure on $X(A)$ such that 
$$\int_{X(A)}f\,d\rho_A=\int_{X(A)}\frac{1}{\#r_A^{-1}(x)}\sum_{r_A(y)=x}f(y)\,d\rho_A(x),$$
for all bounded measurable functions $f$ on $X(A)$.
\par
In this paper we analyze the connection between measures on $X$ 
and the induced measures on $X_\infty$, and we characterize those 
measures $X_\infty$ which are quasi-invariant with respect to the 
invertible mapping $\hat r$
$$\begin{diagram}
X_\infty&\rTo^{\hat r}&X_\infty\\
\dTo& &\dTo\\
X&\rTo^{r}&X
\end{diagram}$$
where
$$X_\infty:=\left\{\hat 
x=(x_0,x_1,\dots )\in\prod_{\bn_0}X\Bigm|r(x_{n+1})=x_n,n\in\bn_0\right\},$$
$$\hat r(\hat x)=(r(x_0),x_0,x_1,\dots )$$
and
$$\hat r^{-1}(\hat x)=(x_1,x_2,\dots ).$$
\par
We need a specific interplay between spaces $(X,r)$, 
$r\colon X\rightarrow X$ a non-invertible endomorphism, and induced 
spaces $(X_\infty,\hat r)$ where $\hat r$ is an automorphism.  
This issue, and variants of it, arise in a number of areas of 
mathematics; first in probability theory, going back to 
\cite{Kol}; and also more recently in a number of wavelet 
problems, see for example \cite{BCMO}, \cite{BM}, \cite{BrJo02}, 
\cite{DutJo}, \cite{DuJo04}, \cite{Gun99}, \cite{Jor01}, and 
\cite{Jor04}. In these applications, the problem is to carry 
along some isometric operator defined on a Hilbert space of 
functions on $X$, to the space $X_\infty$ (see proposition 
\ref{prop2_2} below). In the language of operator theory, we wish 
to make a covariant unitary dilation from $X$ to $X_\infty$. This 
means that we need to induce measures $\mu$ on $X$ to measures 
$\hat\mu$ on $X_\infty$ in such a way that a prescribed 
covariance is preserved.
\par
In our extension of measures from $X$ to $X_\infty$, we must keep 
track of the transfer from one step to the next, and there is an 
operator which accomplishes this, Ruelle's transfer operator (see 
(\ref{eq2_18}) below). In its original form it was introduced in 
\cite{Rue89}, but since that, it has found a variety of 
applications, see e.g., \cite{Bal00}. For use of the Ruelle 
operator in wavelet theory, we refer to \cite{Jor01} and 
\cite{Dut}.

\par
In our construction, the Hilbert spaces of functions on 
$X_\infty$ will be realized as a Hilbert space of martingales. 
This is consistent with our treatment of wavelet resolutions as 
martingales. This was first suggested in \cite{Gun99} in 
connection with wavelet analysis.

\par
To make our paper self-contained, we have recalled Doob's 
martingale convergence theorem in section 3, in the form in which 
we need it, but we refer the reader to the books \cite{Doob3}
and \cite{Neveu} for background on martingale theory.

\section{\label{proj}Projective limits}
\subsection{\label{projdefi}Definitions}
\par
Let $X$ be a compact Hausdorff space, and let $r\colon X\rightarrow X$ 
be a finite-to-one mapping:
\begin{equation}\label{eq2_1}
1\leq\#r^{-1}(x)<\infty
\end{equation}
where $r^{-1}(x)=\{y\in X\,|\,r(y)=x\}$. More generally, set
\begin{equation}\label{eq2_2}
r^{-1}(E)=\{y\in X\,|\,r(y)\in E\},\quad\mbox{if }E\subset X.
\end{equation}
\par
By the projective limit $X_\infty(r)$, we mean
\begin{equation}\label{eq2_3}
X_\infty(r)=\left\{\hat 
x=(x_0,x_1,\dots )\in\prod_{n\in\mathbb{N}_0}X\Bigm|r(x_{n+1})=x_n\,,n\in\mathbb{N}_0\right\}.
\end{equation}
\par
We have the representation
\begin{equation}\label{eq2_4}
\begin{diagram}
X&\lTo^r&X&\lTo^r&X&\lTo^r&\cdots &\lTo X_\infty(r)
\end{diagram}
\end{equation}
\par
We will further assume that $r$ is not invertible, i.e., that 
$\#r^{-1}(x)$ is not the constant function one. \par It is well 
known that, if $r$ is continuous, then pull-backs of open sets in 
$X$ define a topology on $X_\infty(r)$ making $X_\infty(r)$ 
compact.
\par
The restriction to $X_\infty(r)$ of the coordinate projections 
$(x_0,x_1,\dots )\mapsto x_n$ will be denoted by $\theta_n$, and we 
have
\begin{equation}\label{eq2_5}
r\circ\theta_{n+1}=\theta_n\qquad(n\in\mathbb{N}_0).
\end{equation}
\par
One advantage of passing from $X$ to $X_\infty(r)$ is that $r$ 
induces an invertible mapping $\hat 
r\colon X_\infty(r)\rightarrow\xir$, defined by
\begin{equation}\label{eq2_6}
\hat r(\hat x)=\hat r(x_0,x_1,x_2,\dots )=(r(x_0),x_0,x_1,x_2,\dots ).
\end{equation}
One checks that
\begin{equation}\label{eq2_7}
\hat r^{-1}(\hat x)=(x_1,x_2,\dots ),
\end{equation}
i.e., that
\begin{equation}\label{eq2_8}
\hat r\circ\hat r^{-1}=\hat r^{-1}\circ\hat r=\mbox{id}_{\xir},
\end{equation}
and
\begin{equation}\label{eq2_9}
r\circ\theta_n=\theta_n\circ\hat r=\theta_{n-1}.
\end{equation}
\par
Moreover, both $\hat r$ and $\hat r^{-1}$ are continuous if $r$ 
is.
\par
If $\mathfrak{B}$ is a sigma-algebra of subsets in $X$ (typically 
we will take $\mathfrak{B}$ to be the Borel subsets in $X$), then 
there are sigma-algebras
\begin{equation}\label{eq2_10}
\mathfrak{B}_n:=\theta_n^{-1}(\mathfrak{B})=\{\theta_n^{-1}(E)\,|\, 
E\in\mathfrak{B}\}
\end{equation}
and
\begin{equation}\label{eq2_11}
\mathfrak{B}_\infty=\bigcup_{n\in\bn_0}\mathfrak{B}_n.
\end{equation}
These $\mathfrak{B}_n$'s are sigma-algebras of subsets of 
$\xir$.\par Using (\ref{eq2_5}), we get
\begin{equation}\label{eq2_12}
\theta_n^{-1}=\theta_{n+1}^{-1}\circ r^{-1}.
\end{equation}
If $r$ is measurable, then 
$r^{-1}(\mathfrak{B})\subset\mathfrak{B}$, and we conclude that
\begin{equation}\label{eq2_13}
\mathfrak{B}_n\subset\mathfrak{B}_{n+1}.
\end{equation}
Then both of the mappings $\hat r$ and $\hat r^{-1}$ are 
measurable on $\xir$ with respect to the sigma-algebra 
$\mathfrak{B}_\infty$.
\par
In an earlier paper \cite{DuJo04}, we studied measures $\hat\mu$ 
on $(\xir,\mathfrak{B}_\infty)$ and functions
\begin{equation}\label{eq2_14}
V\colon X\rightarrow[0,\infty)
\end{equation}
such that
\begin{equation}\label{eq2_15}
\hat\mu\circ\hat r\ll\hat\mu\qquad(\mbox{quasi-invariance}),
\end{equation}
and
\begin{equation}\label{eq2_16}
\frac{d\hat\mu\circ\hat r}{d\hat\mu}=V\circ\theta_0.
\end{equation}
Notice that when $V$ is given, the measure $\hat\mu$ depends on 
$V$. We say that $\hat\mu$ is {\it $V$-quasi-invariant}.
\par
In \cite{DuJo04} we studied (\ref{eq2_16}) under rather 
restrictive assumptions. Our present study continues the analysis 
of (\ref{eq2_16}), and we give a structure theorem for the 
solutions $(V,\hat\mu)$.
\subsection{\label{projspec}Special solutions to the problem 
(\ref{eq2_15})--(\ref{eq2_16})}
\par
In \cite{DuJo04}, we studied the following restrictive setup: we 
assumed that $X$ carries a probability measure $\mu$ which is 
{\it strongly $r$-invariant}. By this we mean that
\begin{equation}\label{eq2_17}
\int_Xf\,d\mu=\int_X\frac{1}{\#r^{-1}(x)}\sum_{y\in X, 
r(y)=x}f(y)\,d\mu(x)\qquad(f\in L^\infty(X)).
\end{equation}
\par
If, for example $X=\br/\bz$, and $r(x)=2x\mod 1$, then the Haar 
measure on $\br/\bz=$Lebesgue measure on $[0,1)$, is the unique 
strongly $r$-invariant measure on $X$.
\par
Suppose $V$ in (\ref{eq2_14}) is bounded and measurable. Then 
define $R=R_V$, the {\it Ruelle operator}, by
\begin{equation}\label{eq2_18}
R_Vf(x)=\frac{1}{r^{-1}(x)}\sum_{r(y)=x}V(y)f(y)\qquad (f\in 
L^1(X,\mu)).
\end{equation}
\begin{theorem}\label{th2_1}\textup{(\cite{DuJo04})} Let $r\colon X\rightarrow X$ and $\xir$ 
be as described in section \textup{\ref{projdefi}}, and suppose that $X$ 
has a strongly $r$-invariant measure $\mu$. Let $V$ be a 
non-negative, measurable function on $X$, and let $R_V$ be the 
corresponding Ruelle operator.
\begin{enumerate}
\item There is a unique measure $\hat\mu$ on $\xir$ such that
\begin{enumerate}
\item $\hat\mu\circ\theta_0^{-1}\ll\mu$\qquad \textup{(}set 
$h=\frac{d(\hat\mu\circ\theta_0^{-1})}{d\mu}$\textup{)}, \item 
$\displaystyle\int_Xf\,d\hat\mu\circ\theta_n^{-1}=\int_XR^n_V(fh)\,d\mu\qquad(n\in\bn_0).$
\end{enumerate}
\item The measure $\hat\mu$ on $\xir$ satisfies
\begin{equation}\label{eq2_19}
\frac{d(\hat\mu\circ\hat r)}{d\hat\mu}=V\circ\theta_0
\end{equation}
and
\begin{equation}\label{eq2_20}
R_Vh=h.
\end{equation}
\end{enumerate}
\end{theorem}
\par
In this paper we turn around the problem, and take (\ref{eq2_19}) 
as the fundamental axiom. We will not assume the existence of a 
strongly $r$-invariant measure $\mu$.
\begin{proposition}\label{prop2_2}
Let $r\colon X\rightarrow X$ be as described above, and let $\hat 
r\colon\xir\rightarrow\xir$ be the corresponding automorphism. Let 
$m\colon X\rightarrow\bc$ be a bounded measurable function on $X$. Then 
the operator
\begin{equation}\label{eq2_21}
\hat S_mf:=(m\circ\theta_0)f\circ\hat r
\end{equation}
defines an isometry in the Hilbert space $L^2(\xir,\hat\mu)$ if 
and only if
\begin{equation}\label{eq2_22}
\hat\mu\circ\hat r\ll\hat\mu
\end{equation}
and
\begin{equation}\label{eq2_23}
\frac{d(\hat\mu\circ\hat r)}{d\hat\mu}=|m\circ\theta_0|^2.
\end{equation}
\end{proposition}
\begin{proof}
The isometric property for (\ref{eq2_21}) may be stated in the 
form. 
\begin{equation}\label{eq2_24}
\int_{\xir}|m\circ\theta_0|^2|f\circ\hat 
r|^2\,d\hat\mu=\int_{\xir}|f|^2\,d\hat\mu,\mbox{ for all }f\in 
L^2(\xir,\hat\mu).
\end{equation}
Setting $V:=|m|^2$, and $g:=|f\circ\hat r|^2$, (\ref{eq2_24}) 
reads
$$\int_{\xir}V\circ\theta_0\,g\,d\hat\mu=\int_{\xir}g\circ\hat 
r\,d\hat\mu.$$ But this amounts precisely to the two assertions 
(\ref{eq2_22}) and (\ref{eq2_23}). In other words, the 
$V$-quasi-invariance property for $V=|m|^2$ is equivalent to 
$S_m$ defining an $L^2$-isometry.
\end{proof}
\section{\label{indu}Induction of measures on $\xir$ from measures 
on $X$}
\subsection{\label{indumart}Martingales}
\par
If $\hat\mu$ is a measure defined on the sigma-algebra 
$\mathfrak{B}_\infty$ on $\xir$, then there is an associated 
sequence of measures $(\mu_n)$ on $X$ defined as follows:
\begin{equation}\label{eq3_1}
\mu_n:=\hat\mu\circ\theta_n^{-1}\quad(n\in\bn_0),
\end{equation}
or more precisely
\begin{equation}\label{eq3_2}
\mu_n(E)=\hat\mu(\theta_n^{-1}(E))\qquad(E\in\mathfrak{B}),
\end{equation}
where
\begin{equation}\label{eq3_3}
\theta_n^{-1}(E)=\{\hat x\in\xir\,|\,\theta_n(\hat x)\in E\}.
\end{equation}
Our measures will be assumed positive and finite, unless 
specified otherwise.
\par
We now introduce the Hilbert spaces:
\begin{equation}\label{eq3_4}
\hat H(\hat\mu):=L^2(\xir,\hat\mu)
\end{equation}
and
\begin{equation}\label{eq3_5}
H_n(\hat\mu):=\{\xi\circ\theta_n\,|\,\xi\in L^2(X,\mu_n)\},
\end{equation}
the orthogonal projections
\begin{equation}\label{eq3_6}
P_n\colon\hat H(\hat\mu)\rightarrow H_n(\hat\mu)
\end{equation}
and
\begin{equation}\label{eq3_7}
E_n\colon\hat H(\hat\mu)\rightarrow L^2(X,\mu_n),
\end{equation}
setting
\begin{equation}\label{eq3_8}
P_n(f)=:E_n(f)\circ\theta_n\qquad(n\in\bn_0);
\end{equation}
then
\begin{equation}\label{eq3_9}
\int_{\xir}\xi\circ\theta_n\,f\,d\hat\mu=\int_X\xi\, 
E_n(f)\,d\mu_n,\mbox{ for all }n\in\bn_0, \xi\in L^2(X,\mu_n).
\end{equation}
\begin{lemma}\label{lem3_1}
Let $\hat\mu$ be a measure on $(\xir,\mathfrak{B}_\infty)$ and 
let $\hat H(\hat\mu)$be the corresponding Hilbert space from 
\textup{(\ref{eq3_4})}.
\begin{enumerate}
\item Then there is an isometric isomorphism 
\begin{equation}\label{eq3_10}
J\colon\hat H(\hat\mu)\rightarrow\lim_n L^2(X,\mu_n)
\end{equation}
where the inductive limit in \textup{(\ref{eq3_10})} is defined by the 
isometric embeddings
\begin{equation}\label{eq3_11}
\xi\mapsto\xi\circ r\,\colon L^2(X,\mu_n)\rightarrow L^2(X,\mu_{n+1}).
\end{equation}
The isomorphism $J$ in \textup{(\ref{eq3_10})} is
\begin{equation}\label{eq3_12}
Jf=(E_nf)_{n\in\bn_0}\qquad(f\in\hat H(\hat\mu)).
\end{equation}
\item The Hilbert norm is given by $$\|f\|_{\hat 
H(\hat\mu)}^2=\lim_{n\rightarrow\infty}\|E_n(f)\|_{L^2(X,\mu_n)}^2;$$
and \item For every sequence $\xi_n\in L^2(X,\mu_n)$ such that 
$P_n(\xi_{n+1}\circ\theta_{n+1})=\xi_n\circ\theta_n$, and
\begin{equation}\label{eq3_13}
\sup_n\|\xi_n\|_{L^2(X,\mu_n)}<\infty,
\end{equation}
there is a unique $f\in L^2(\xir,\hat\mu)$ such that
\begin{equation}\label{eq3_14}
P_nf=\xi_n\circ\theta_n\mbox{ for all }n\in\bn_0.
\end{equation}
Moreover,
\begin{equation}\label{eq3_15}
f=\lim_{n\rightarrow\infty}\xi_n\circ\theta_n,\mbox{ pointwise 
}\hat\mu\mbox{-a.e., and in }\hat H(\hat\mu).
\end{equation}
\end{enumerate}
\end{lemma}
\begin{proof}
The proof depends on Doob's martingale convergence theorem, see 
\cite{DuJo04}, \cite{Neveu}, and \cite{Doob3}.
\end{proof}
\subsection{\label{indufrom}From $\mu$ on $X$ to $\hat\mu$ on 
$X_\infty(r)$}
\begin{lemma}\label{lem3_2}
Let $(\mu_n)_{n\in\bn_0}$ be a sequence of measures on 
$(X,\mathfrak{B})$. Then there is a measure $\hat\mu$ on 
$(\xir,\mathfrak{B}_\infty)$ such that 
\begin{equation}\label{eq3_16}
\hat\mu\circ\theta_n^{-1}=\mu_n\qquad(n\in\bn_0)
\end{equation}
if and only if
\begin{equation}\label{eq3_17}
\mu_{n+1}\circ r^{-1}=\mu_n\qquad(n\in\bn_0).
\end{equation}
\end{lemma}
\begin{proof}
We first recall the definition, the measure $\mu_{n+1}\circ 
r^{-1}$ in (\ref{eq3_17}). On a Borel subset $E$, it is 
\begin{equation}\label{eq3_17new}
\mu_{n+1}\circ r^{-1}(E):=\mu_{n+1}(r^{-1}(E)).
\end{equation}
It is well known (see \cite{Rud87}) that (\ref{eq3_17new}) 
defines a measure.
\par
 Suppose $\hat\mu$ exists such that (\ref{eq3_16}) holds. Let 
$\xi$ be a bounded $\mathfrak{B}$-measurable function on $X$. Then
$$\int_X\xi\circ r\,d\mu_{n+1}=\int_{\xir}\xi\circ 
r\circ\theta_{n+1}\,d\hat\mu=\int_{\xir}\xi\circ\theta_n\,d\hat\mu=\int_X\xi\,d\hat\mu\circ\theta_n^{-1}=\int_X\xi\,d\mu_n,$$
which proves the identity (\ref{eq3_17}).
\par
Conversely, suppose $(\mu_n)_{n\in\bn_0}$ satisfies 
(\ref{eq3_17}). Note that the $\mathfrak{B}_n$-measurable 
functions are of the form $f\circ\theta_n$ with $f$ 
$\mathfrak{B}$-measurable. The $\mathfrak{B}_n$-measurable 
functions are also $\mathfrak{B}_{n+1}$-measurable and this 
inclusion is given by 
\begin{equation}\label{eq3_18}
f\circ\theta_n=(f\circ r)\circ\theta_{n+1}.
\end{equation}
\par
The relations (\ref{eq2_11}) and (\ref{eq2_13}) imply that the 
union of the algebras of bounded $\mathfrak{B}_n$-measurable 
functions is $L^1(\hat\mu)$-dense in the algebra of bounded 
$\mathfrak{B}_\infty$-measurable functions on $\xir$.
\par
Now set
\begin{equation}\label{eq3_19}
\hat\mu(f\circ\theta_n):=\mu_n(f),\, f\mbox{ assumed 
}\mathfrak{B}_n\mbox{-measurable and bounded.}
\end{equation}
To see that this is consistent, use (\ref{eq3_18}) and 
(\ref{eq3_17}) to check that
$$\hat\mu((f\circ r)\circ\theta_{n+1})=\mu_{n+1}(f\circ 
r)=\mu_n(f)=\hat\mu(f\circ\theta_n).$$ The existence and 
uniqueness of $\hat\mu$ follows from Kolmogorov's theorem, see 
\cite{DuJo04} and \cite{Kol}.
\end{proof}
\subsection{\label{indumain}The main theorem}
\par
In this section we will study the problem (\ref{eq2_16}) via our 
correspondence $\hat\mu\leftrightarrow(\mu_n)_{n\in\bn_0}$ from 
lemma \ref{lem3_2}.
\par
The questions are two: (1) Which measures $\mu_0$ on $X$ admit 
``extensions'' $\hat\mu$ to measures on $\xir$ which satisfy 
(\ref{eq2_16})?
\par
(2) Understand the measures $\mu_0$ which admit solutions 
$\hat\mu$ to (\ref{eq2_16}).
\par
Our first theorem answers the question (1).
\begin{theorem}\label{th3_3}
Let $(\mu_n)_{n\in\bn_0}$ be a sequence of measures as in lemma 
\textup{\ref{lem3_1}}, i.e., satisfying $\mu_{n+1}\circ r^{-1}=\mu_n$, and 
let $V\colon X\rightarrow[0,\infty)$ be $\mathfrak{B}$-measurable. Then 
the ``extended'' measure $\hat\mu$ satisfies
\begin{equation}\label{eq3_20}
\frac{d(\hat\mu\circ\hat r)}{d\hat\mu}=V\circ\theta_0\mbox{ a.e. 
on }(\xir,\mathfrak{B}_\infty)
\end{equation}
if and only if
\begin{equation}\label{eq3_21}
d\mu_0=(V\,d\mu_0)\circ r^{-1},\mbox{ and }\,d\mu_{n+1}=(V\circ 
r^n)d\mu_n\qquad(n\in\bn_0),
\end{equation}
or equivalently
\begin{equation}\label{eq3_22}
\!\int_Xf\,d\mu_0=\int_Xf\circ r\,V\,d\mu_0,\mbox{ and}
\int_Xf\,d\mu_{n+1}=\int_Xf\,V\circ r^n\,d\mu_n\qquad(n\in\bn_0)
\end{equation}
for all bounded $\mathfrak{B}$-measurable functions $f$ on $X$.
\end{theorem}

\section{\label{proof}Proof of theorem \ref{th3_3}}
Let $(\mu_n)_{n\in\bn_0}$ be a sequence of measures that satisfy
\begin{equation}\label{eq4_1}
\mu_{n+1}\circ r^{-1}=\mu_n,
\end{equation}
and let $\hat\mu$ be the corresponding measure on $\xir$. Suppose 
first that $\hat\mu$ satisfies (\ref{eq3_20}) for some 
$\mathfrak{B}$-measurable function $V\colon X\rightarrow[0,\infty)$; in 
other words we have the identity
\begin{equation}\label{eq4_2}
\int_{\xir} f\circ\hat 
r^{-1}\,d\hat\mu=\int_{\xir}f\,V\circ\theta_0\,d\hat\mu.
\end{equation}
for all bounded $\mathfrak{B}_\infty$-measurable functions $f$ on 
$\xir$.
\par
Now specializing (\ref{eq4_2}) to $f=\xi\circ\theta_n$ with $\xi$ 
some bounded $\mathfrak{B}$-measurable function on $X$, we get
\begin{align*}
\int_{\xir}\xi\circ\theta_{n+1}\,d\hat\mu&=\int_{\xir}\xi\circ\theta_n\,V\circ\theta_0\,d\hat\mu\\
&=\int_{\xir}\xi\circ\theta_n\,V\circ 
r^n\circ\theta_n\,d\hat\mu=\int_X\xi\,V\circ r^n\,d\mu_n,
\end{align*}
and 
therefore
\begin{equation}\label{eq4_3}
d\mu_{n+1}=(V\circ r^n)\,d\mu_n
\end{equation}
by (\ref{eq3_1}).
\par
We take a closer look at the case $n=0$. Substitute 
$f=\xi\circ\theta_0$ into (\ref{eq4_2}). We get
$$\int_{\xir}\xi\circ\theta_0\,d\hat\mu=\int_{\xir}\xi\circ\theta_0\circ\hat 
r\,V\circ\theta_0\,d\hat\mu=\int_{\xir}\xi\circ 
r\circ\theta_0\,V\circ\theta_0\,d\hat\mu=\int_X\xi\circ 
r\,V\,d\mu_0$$ and therefore $d\mu_0=(V\,d\mu_0)\circ r^{-1}$ 
which is (\ref{eq3_21}).
\par
Conversely, suppose $(\mu_n)_{n\in\bn_0}$ satisfies (\ref{eq4_1}) 
and (\ref{eq3_21}). Lemma \ref{lem3_2} yields the existence of a 
$\hat\mu$ on $(\xir,\mathfrak{B}_\infty)$ such that 
\begin{equation}\label{eq4_4}
\mu_n=\hat\mu\circ\theta_n^{-1}.
\end{equation}
\par
We claim that $\hat\mu$ is $V$-quasi-invariant, i.e., that 
(\ref{eq3_20}) holds. This amounts to the identity
\begin{equation}\label{eq4_5}
\int_{\xir}f\circ\hat 
r^{-1}\,d\hat\mu=\int_{\xir}f\,V\circ\theta_0\,d\hat\mu.
\end{equation}
By lemma \ref{lem3_1}, this is equivalent to 
$$d\mu_{n+1}=V\circ r^n\,d\mu_n.$$
\par
To see this, substitute $f=\xi\circ\theta_n$ into (\ref{eq4_5}) 
for $\xi$ a bounded $\mathfrak{B}$-measurable function on $X$, 
and $n\in\bn_0$. But this holds by (\ref{eq3_21}), and 
(\ref{eq3_20}) follows.
\section{\label{fixe}A fixed-point problem}
\par
In theorem \ref{th3_3}, we saw that our condition (\ref{eq2_16}), 
{\it $V$-quasi-invariance}, on measures $\hat\mu$ on 
$(\xir,\mathfrak{B}_\infty)$ entails a {\it fixed-point property} 
(\ref{eq5_no}) for the corresponding measure 
$\mu_0:=\hat\mu\circ\theta_0^{-1}$. In this section, we turn the 
problem around. We show that this fixed-point property 
characterizes the measures $\hat\mu$ satisfying the 
quasi-invariance, i.e., satisfying (\ref{eq2_16}).
\par
If the function $V\colon X\rightarrow[0,\infty)$ is given, we define 
$$V^{(n)}(x):=V(x)V(r(x))\cdots V(r^{n-1}(x)),$$
and set $d\mu_n:=V^{(n)}d\mu_0$. Our result states that the 
corresponding measure $\hat\mu$ on $\xir$ is $V$-quasi-invariant 
if and only if 
\begin{equation}\label{eq5_no}
d\mu_0=(V\,d\mu_0)\circ r^{-1}.
\end{equation}

\begin{theorem}\label{th5_1}
Let $V\colon X\rightarrow[0,\infty)$ be $\mathfrak{B}$-measurable, and 
let $\mu_0$ be a measure on $X$ satisfying the following 
fixed-point property
\begin{equation}\label{eq5_1}
d\mu_0=(V\,d\mu_0)\circ r^{-1}.
\end{equation}
Then there exists a unique measure $\hat\mu$ on $\xir$ such that 
\begin{equation}\label{eq5_2}
\frac{d(\hat\mu\circ\hat r)}{d\hat\mu}=V\circ\theta_0
\end{equation}
and
$$\hat\mu\circ\theta_0^{-1}=\mu_0.$$
\end{theorem}

\begin{proof}
Let $\mu_0$ be a measure satisfying the fixed-point property 
(\ref{eq5_1}), and set
\begin{equation}\label{eq5_3}
d\mu_{n+1}:=(V\circ r^n)\,d\mu_n.
\end{equation}
This defines a sequence $\mu_0,\mu_1,\mu_2,\dots $ of measures on 
$X$. To establish the existence of the measure $\hat\mu$ on 
$X_\infty(r)$, we appeal to lemma \ref{lem3_2}. We claim that 
(\ref{eq3_17}) holds, or equivalently that
\begin{equation}\label{eq5_4}
\int_Xf\circ r\,d\mu_{n+1}=\int_Xf\,d\mu_n
\end{equation}
for all bounded $\mathfrak{B}$-measurable functions $f$ on $X$, 
and for all $n\in\bn_0$.
\par
Substitution of (\ref{eq5_3}) into (\ref{eq5_4}) yields
$$\int_Xf\circ r\,d\mu_{n+1}=\int_Xf\circ r\,V\circ r^n\,d\mu_n,$$
and we will prove (\ref{eq5_4}) by induction:
\par
First
$$\int_Xf\circ r\,d\mu_1=\int_Xf\circ 
r\,V\,d\mu_0=\int_Xf\,(V\,d\mu_0)\circ r^{-1}=\int_Xf\,d\mu_0,$$ 
where we used (\ref{eq5_1}) in the last step. So (\ref{eq5_4}) 
holds for $n=0$.
\par
Assume that (\ref{eq5_4}) holds for $\mu_m$ and $m<n$. Then 
$$\int_Xf\circ r\,d\mu_{n+1}=\int_Xf\circ r\,V\circ r^{n-1}\circ 
r\,d\mu_n=\int_Xf\,V\circ r^{n-1}\,d\mu_{n-1}=\int_Xf\,d\mu_n,$$ 
where we used the definition (\ref{eq5_3}) in the last step. The 
induction is completed.
\par
This proves (\ref{eq3_17}). An application of lemma \ref{lem3_2} 
yields the existence of a unique measure $\hat\mu$ on 
$(\xir,\mathfrak{B}_\infty)$ such that 
\begin{equation}\label{eq5_5}
\hat\mu\circ\theta_n^{-1}=\mu_n\qquad(n\in\bn_0).
\end{equation}
Since the sequence $(\mu_n)_{n\in\bn_0}$ satisfies the pair of 
conditions (\ref{eq3_21}) by construction, we conclude from 
theorem \ref{th3_3} that $\hat\mu$ must satisfy (\ref{eq5_2}), 
i.e., $\hat\mu$ is $V$-quasi-invariant as claimed. 
\end{proof}
\par
The proof of theorem \ref{th5_1} yields also the following lemma
\begin{lemma}\label{lemno_1}
Let $V\colon X\rightarrow[0,\infty)$ be a bounded 
$\mathfrak{B}$-measurable function, and let $\mu_0$ be a measure 
on $(X,\mathfrak{B})$. Suppose
\begin{equation}\label{eqno_1_1}
d\mu_0=(V\,d\mu_0)\circ r^{-1}.
\end{equation}
Set
\begin{equation}\label{eqno_1_3}
d\mu_n:=V^{(n)}\,d\mu_0\mbox{ for }n\in\bn.
\end{equation}
then
\begin{equation}\label{eqno_1_4}
\mu_{n+1}\circ r^{-1}=\mu_n,
\end{equation}
and
\begin{equation}\label{eqno_1_5}
\mu_n(X)=\int_XV^{(n)}\,d\mu_0=\mu_0(X),\mbox{ for all }n\in\bn.
\end{equation}
\end{lemma}
\begin{proof}
Indeed, by (\ref{eqno_1_3}) we have the recursive formula 
$d\mu_{n+1}=(V\circ r^n)\,d\mu_n$, and as in the proof of theorem 
\ref{th5_1}, we get that (\ref{eqno_1_4}) is satisfied. Applying 
(\ref{eqno_1_3}) to the constant function $1$, we get that 
$$\mu_n(X)=\int_XV^{(n)}\,d\mu_0,$$
and, with (\ref{eqno_1_4}) 
$$\mu_{n+1}(X)=\int_X1\,d\mu_{n+1}=\int_X1\circ 
r\,d\mu_{n+1}=\int_X1\,d\mu_n=\mu_n(X).$$ Now (\ref{eqno_1_5}) 
follows by induction.
\end{proof}
\section{\label{tran}Transformations of measures}
Let $X$ be a compact Hausdorff space, $\mathfrak{B}$ the 
sigma-algebra of all Borel subsets of $X$, $C(X)$ the continuous 
functions on $X$. 
\par
For $a\in\br_+$ denote by
$$M_a(X):=\{\mu\,|\,\mu\mbox{ is a measure on 
}(X,\mathfrak{B}),\mu(X)=a\}.$$ Note that $M_a(X)$ is equipped 
with the weak$^*$-topology coming from the duality
$$(f,\mu)\mapsto\int_Xf\,d\mu=:\mu(f)$$
where $f\in C(X)$. The neighborhoods are generated by the sets 
$N_{\mu_0}(f_1,\dots ,f_k,\epsilon)$ where $\epsilon\in\br_+$, 
$f_1,\dots ,f_k\in C(X)$ and 
$$N_{\mu_0}(f_1,\dots ,f_k,\epsilon)=\{\mu\in M(X)\,|\, 
|\mu(f_i)-\mu_0(f_i)|<\epsilon, i=1,\dots,k\}$$ \par It is known 
\cite{Rud87} that each $M_a(X)$ is weak$^*$-compact, i.e., 
$M_a(X)$ is a compact convex set in the topology determined by 
the neighborhoods $N_{\mu_0}(f_1,\dots ,f_k,\epsilon)$.
\par
Let $V\colon X\rightarrow[0,\infty)$ be bounded and 
$\mathfrak{B}$-measurable, and define
$$T_V(\mu):=(V\,d\mu)\circ r^{-1}\qquad(\mu\in M(X)).$$
\begin{lemma}\label{lem6_1}
Let $V$ be as above, and assume $V$ is also continuous. Then 
\linebreak 
$T_V\colon M(X)\rightarrow M(X)$ is continuous. 
\end{lemma}
\begin{proof}
Let $\mu\in M(X)$, $k\in\bz_+$ and $f_i\in C(X)$, $i=1,\dots ,k$. 
Then 
$$\ip{T_V(\mu)}{f_i}=\ip{\mu}{V\,f_i\circ r}.$$
Setting $g_k:=V\,f_i\circ r$, we see that
$$T_V(N_{\mu_0}(g_1,\dots ,g_k,\epsilon))\subset 
N_{T_V(\mu_0)}(f_1,\dots ,f_k,\epsilon)$$ and the conclusion follows 
if $V$ is assumed continuous.
\end{proof}
\begin{definition}\label{def6_no}
Let $V\colon X\rightarrow[0,\infty)$ be bounded and 
$\mathfrak{B}$-measurable. We use the notation
$$M^V(X):=\{\mu\in M(X)\,|\,d\mu=(V\,d\mu)\circ r^{-1}\}.$$
\par
For measures $\hat\mu$ on $(\xir,\mathfrak{B}_\infty)$ we 
introduce 
$$M_{qi}^V(\xir):=\left\{\hat\mu\in M(\xir)\Bigm|\hat\mu\circ\hat 
r\ll\hat\mu\mbox{ and }\frac{d(\hat\mu\circ\hat 
r)}{d\hat\mu}=V\circ\theta_0\right\}.$$
\end{definition}
\par
The results of the previous section may be summarized as follows:
\begin{theorem}\label{th6_no_no}
Let $V$ be as in definition \textup{\ref{def6_no}}. For measures $\hat\mu$ 
on $\xir$ and $n\in\bn_0$, define 
$$C_n(\hat\mu):=\hat\mu\circ\theta_n^{-1}.$$
then $C_0$ is a bijective affine isomorphism of $M_{qi}^V(\xir)$ 
onto $M^V(X)$ that preserves the total measure, i.e., 
$C_0(\hat\mu)(X)=\hat\mu(\xir)$ for all $\hat\mu\in 
M_{qi}^V(\xir)$.
\end{theorem}
\begin{proof}
We showed in theorem \ref{th3_3} that $C_0$ maps $M_{qi}^V(\xir)$ 
onto $M^V(X)$. The inverse mapping 
$$C_0^{-1}\colon M^V(X)\rightarrow M_{qi}^V(\xir)$$
may be realized using theorem \ref{th5_1} and lemma \ref{lemno_1}.
\end{proof}
\begin{remark}\label{rem6_3}
The intuitive idea behind $C_0^{-1}\colon M^V(X)\rightarrow 
M_{qi}^V(\xir)$ is as follows: 
\par
Let $C_0^{-1}(\mu_0)=\hat\mu.$ Then $\hat\mu$ is a measure on 
$\xir$, and we view points $\hat x$ in $\xir$ as infinite paths. 
Recall, if $\hat x=(x_0,x_1,\dots )\in\xir$ then $r(x_{n+1})=x_n$. 
So in a random walk we choose $x_{n+1}\in r^{-1}(x_n)$, and the 
function $V$ assigns the probabilities in each step. At step $n$, 
there are $\#r^{-1}(x_n)$ choices.
\par
The assertion in the theorem is that the measure $\hat\mu$ is 
completely determined by the prescribed measure $\mu_0$ at the 
starting point $x_0$ of the path, and by the function $V$.
\par
The two measures $\hat\mu$ and $\mu_0$ are normalized so that 
$$\mu_0(X)=\hat\mu(\xir).$$
\par
In a special case, we give an explicit formula for 
$C_0^{-1}(\mu_0)$, \cite[Proposition 8.2]{DuJo04}
\end{remark}
\begin{theorem}\label{th6_nooo}
Let $V\colon X\rightarrow[0,\infty)$ be continuous. Assume that there 
exist some measure $\nu$ on $(X,\mathfrak{B})$ and two numbers 
$0<a<b$ such that 
\begin{equation}\label{eq6_nooo_1}
a\leq\nu(X)\leq b, \mbox{ and }a\leq\int_X V^{(n)}\,d\nu\leq 
b\mbox{ for all }n\in\bn.
\end{equation}
Then there exists a measure $\mu_0$ on $(X,\mathfrak{B})$ that 
satisfies 
$$d\mu_0=(V\,d\mu_0)\circ r^{-1},$$
and there exists a $V$-quasi-invariant measure $\hat\mu$ on 
$(\xir,\mathfrak{B}_\infty)$.
\end{theorem}
\begin{proof}
Condition (\ref{eq6_nooo_1}) guarantees that the set 
$$M_{ab}^V(X):=\left\{\mu\in M(X)\Bigm|a\leq\mu(X)\leq b\mbox{ and 
}a\leq\int_XV^{(n)}\,d\mu\leq b,\mbox{ for all }n\in\bn\right\}$$ is 
non-empty. Moreover, by the Banach-Alaoglu theorem \cite{Rud87}, 
this set is compact in the weak$^*$-topology. It is clear also 
that the set is convex.
\par
We claim that the operator $T_V$ maps $M_{ab}^V(X)$ into itself. 
Indeed, if $\mu\in M_{ab}^V(X)$, then 
$$T_V(\mu)(X)=\int_XV\,1\circ r\,d\mu\in[a,b]$$
and, using  $V\,V^{(n)}\circ r=V^{(n+1)}$,
$$T_V(\mu)(V^{(n)})=\int_XV\,V^{(n)}\circ r\,d\mu=\int_XV^{(n+1)}\,d\mu\in[a,b].$$
\par
Lemma \ref{lem6_1} shows that $T_V$ is continuous so we can apply 
the Markov-Kakutani  fixed-point theorem (see \cite{Rud91}) to 
conclude that there exists a $\mu_0\in M_{ab}^V$ such that 
$T_V(\mu_0)=\mu_0$. The $V$-quasi-invariant measure $\hat\mu$ can 
be obtained from $\mu_0$ using theorem \ref{th5_1}.
\end{proof}

\begin{remark}\label{rem6_5}
We now connect the above discussion with section \ref{projspec}. Let 
$(X,\mathfrak{B})$ and $r\colon X\rightarrow X$ be as described above. 
Suppose in addition that $(X,\mathfrak{B})$ carries a strongly 
$r$-invariant probability measure, $\rho$, see (\ref{eq2_17}) for 
the definition.
\par
Set $$M(\rho,X):=\{\mu\in M(X)\,|\,\mu\ll\rho\}.$$ Let 
$V\colon X\rightarrow [0,\infty)$ be a bounded measurable function. 
Then $T_V$ leaves $M(\rho,X)$ invariant; and  $\mu_0$ in 
$M(\rho,X)$ satisfies
$$T_V(\mu_0)=\mu_0$$
if and only if there is an $h\in L^1(X,\rho)$ such that 
$d\mu_0=h\,d\rho$ and $R_Vh=h$.
\end{remark}
\begin{proof}
Let $\mu\in M(\rho,X)$, and write $d\mu=f\,d\rho$, $f\in 
L^1(X,\rho)$. Let $\xi$ be a bounded measurable function on $X$. 
Then
$$T_V(\mu)(\xi)=\int_X\xi\circ 
r\,V\,f\,d\rho=\int_X\xi(x)\frac{1}{\#r^{-1}(x)}\sum_{r(y)=x}V(y)f(y)\,d\rho(x)$$$$=\int_X\xi(x)(R_Vf)(x)\,d\rho(x).$$
Stated differently,
$$T_V(f\,d\rho)=R_V(f)\,d\rho.$$
Hence $T_V(f\,d\rho)=f\,d\rho$ if and only if $R_Vf=f$ as claimed.
\end{proof}
\par
Before stating and proving our next result we need a lemma.
\begin{lemma}\label{lem6_6}
Let $(X,\mathfrak{B})$, $r\colon X\rightarrow X$, and 
$V\colon X\rightarrow[0,\infty)$ be as above. Suppose 
$(X,\mathfrak{B})$ carries a strongly $r$-invariant probability 
measure $\rho$. Then 
\begin{equation}\label{eq6_6_1}
\int_XV^{(n)}\,f\,d\rho=\int_XR_V^n(f)\,d\rho
\end{equation}
for all bounded measurable functions $f$ and $n\in\bn$.
\end{lemma}
\begin{proof}
We prove (\ref{eq6_6_1}) by induction, starting with $n=1$. 
Indeed, for $f\in L^\infty(X)$,
$$\int_XV\,f\,d\rho=\int_X\frac{1}{\#r^{-1}(x)}\sum_{r(y)=x}V(y)f(y)\,d\rho(x)=\int_XR_Vf\,d\rho$$
where we used the definition (\ref{eq2_17}) of strong 
$r$-invariance. 
\par
Suppose (\ref{eq6_6_1}) holds up to $n$. Then 
\begin{align*}
\int_XV^{(n+1)}\,f\,d\rho&=\int_X(V^{(n)}\circ r)Vf\,d\rho=\int_X 
V^{(n)}(x)\frac{1}{\#r^{-1}(x)}\sum_{r(y)=x}V(y)f(y)\,d\rho(x)\\
&=\int_XV^{(n)}R_Vf\,d\rho=\int_XR_V^nR_Vf\,d\rho=\int_XR_V^{n+1}f\,d\rho,
\end{align*}
where we used the induction hypothesis in the last step.
\par
This completes the proof of the lemma.
\end{proof}
\begin{theorem}\label{th6_7}
Let $(X,\mathfrak{B})$, and $r\colon X\rightarrow X$, be as described 
above. Suppose $V\colon X\rightarrow[0,\infty)$ is measurable, 
$$\frac{1}{\#r^{-1}(x)}\sum_{r(y)=x}V(y)\leq1,$$
and that some probability measure $\nu_V$ on $X$ satisfies 
\begin{equation}\label{eq6_7_no}
\nu_V\circ R_V=\nu_V.
\end{equation} Assume also that 
$(X,\mathfrak{B})$ carries a strongly $r$-invariant probability 
measure $\rho$, such that
\begin{equation}\label{eq6_7_2}
\rho(\{x\in X\,|\, V(x)>0\})>0.
\end{equation}
Then
\begin{enumerate}
\item $T_V^n(\,d\rho)=R_V^n(\mathbf{1})\,d\rho,$ for $n\in\bn,$ 
where $\mathbf{1}$ denotes the constant function one. \item 
 {\textup{[}Monotonicity\/\textup{]}} $\quad \cdots \leq R_V^{n+1}(\mathbf{1})\leq 
R_V^n(\mathbf{1})\leq\dots \leq\mathbf{1}$, pointwise on $X$.

\item The limit $\lim_{n\rightarrow}R_V^n(\mathbf{1})=h_V$ 
exists, $R_Vh_V=h_V$, and 
\begin{equation}\label{eq6_7_3}
\rho(\{x\in X\,|\, h_V(x)>0\})>0; 
\end{equation}
\item The measure $d\mu_0^{(V)}=h_V\,d\rho$ is a solution to the 
fixed-point problem $$T_V(\mu_0^{(V)})=\mu_0^{(V)}.$$ \item The 
sequence $d\mu_n^{(V)}=V^{(n)}\,h_V\,d\rho$ defines a unique 
$\hat\mu^{(V)}$ as in theorem \textup{\ref{th5_1}} and lemma 
\textup{\ref{lemno_1}}; and \item 
$\mu_n^{(V)}(f)=\int_XR_V^n(fh_V)\,d\rho$ for all bounded 
measurable functions $f$ on $X$, and all $n\in\bn$.
\par
Finally, \item The measure $\hat\mu^{(V)}$ on $\xir$ satisfying 
$\hat\mu^{(V)}\circ\theta_n^{-1}=\mu_n^{(V)}$ has total mass
$$\hat\mu^{(V)}(\xir)=\rho(h_V)=\int_Xh_V(x)\,d\rho(x).$$
\end{enumerate}
\end{theorem}
\begin{proof}
Part (i) follows from remark \ref{rem6_5}. 
\par
(ii) It is clear that
$$R_V(\mathbf{1})=\frac{1}{\#r^{-1}(x)}\sum_{r(y)=x}V(y)\leq1,\mbox{ 
for }x\in X.$$ Starting the induction, suppose
$$R_V^n(\mathbf{1})\leq R_V^{n-1}(\mathbf{1}).$$
Then \begin{align*}
R_V^{n+1}(\mathbf{1})&=\frac{1}{\#r^{-1}(x)}\sum_{r(y)=x}V(y)R_V^n(\mathbf{1})(y)\\
&\leq\frac{1}{\#r^{-1}(x)}\sum_{r(y)=x}V(y)R_V^{n-1}(\mathbf{1})(y)=R_V^n(\mathbf{1})(x).
\end{align*}
An induction now proves (ii).
\par
It follows that the limit $h_V$ in (iii) exists, and that 
\begin{equation}\label{eq6_7_4}
R_V(h_V)=h_V.
\end{equation}
Using (\ref{eq6_7_no}), we get that 
$$\nu_V(R_V^n(\mathbf{1}))=\nu_V(\mathbf{1})=1,\mbox{ for all 
}n\in\bn,$$ and therefore 
$$\nu_V(h_V)=\int_Xh_V\,d\nu_V=1.$$
 The conclusion (\ref{eq6_7_3}) follows from \cite[Chapter 
3]{Jor04}. 
\par
The assertion in (iv) is that $\mu_0^{(V)}:=h_V\,d\rho$ is a 
fixed point, i.e., that
$$T_V(h_V\,d\rho)=h_V\,d\rho.$$
This follows in turn from remark \ref{rem6_5}. Hence, the 
measures $$d\mu_n^{(V)}=V^{(n)}h_V\,d\rho$$ extend to a unique 
measure $\hat\mu^{(V)}$ on $\xir$ and, by theorem 
\ref{th6_no_no},$$\hat\mu^{(V)}\in M_{qi}^V(\xir).$$ 
\par
Moreover, using now lemma \ref{lem6_6}, we get 
$$\mu_n^{(V)}(f)=\int_Xf\,V^{(n)}\,h_V\,d\rho=\int_XR_V^n(fh_V)\,d\rho,$$
which yields the desired conclusion (vi) and (vii). \par Indeed
$$\hat\mu^{(V)}(\xir)=\int_{\xir}\mathbf{1}\circ\theta_0\,d\hat\mu^{(V)}=\int_X\mathbf{1}\,d\mu_0=\int_Xh_V\,d\rho.$$
Comparing with (\ref{eq6_nooo_1}) in theorem \ref{th6_nooo}, 
notice that
$$\int_Xh_V\,d\rho=\lim_{n\rightarrow\infty}\int_XR_V^n(\mathbf{1})\,d\rho=\lim_{n\rightarrow\infty}\rho(V^{(n)}).$$
\end{proof}
The conclusions of theorem \ref{th6_7} hold in the following more general 
form.
\begin{corollary}\label{cor6_8}
Let $V,W\colon X\rightarrow[0,\infty)$ be measurable, and suppose the 
assumptions \textup{(i)} and \textup{(ii)} are satisfied:
\begin{enumerate}
\item $\displaystyle\frac{1}{\#r^{-1}(x)}\sum_{r(y)=x}V(y)W(y)\leq1$, $x\in X$.
\item
 There exists a non-negative function $h_W$ on $X$, and a measure 
 $\rho_W$ such that $\rho_WR_W=\rho_W$, $R_Wh_W=h_W$ and 
 $\rho_W(h_W)=1$.
 \end{enumerate}
 Then
 \begin{equation}\label{eq6_8_1}
T_V^n(f\,d\rho_W)=(R_{VW}^nf)\,d\rho_W,
\end{equation}
 and the limit 
 $$\lim_{n\rightarrow\infty}R_{VW}^n(\mathbf{1})=h$$
 exists, and satisfies
  \begin{equation}\label{eq6_8_2}
T_V(h\,d\rho_W)=h\,d\rho_W.
\end{equation}
\end{corollary}
\begin{proof}
We first show (\ref{eq6_8_1}) for $n=1$. Let $f$ and $\xi$ be 
bounded, measurable functions on $X$. Then
$$T_V(f\,d\rho_W)(\xi)=\int_X\xi\circ 
r\,V\,f\,d\rho_W=\int_XR_W(\xi\circ 
r\,V\,f)\,d\rho_W$$$$=\int_X\xi\,R_W(Vf)\,d\rho_W=\int_X\xi\,R_{VW}(f)\,d\rho_W.$$

Hence (\ref{eq6_8_1}) holds for $n=1$, and the general case 
follows by iteration.
\par
The argument from the proof of the theorem shows that
$$\cdots \leq R_{VW}^{n+1}(\mathbf{1})\leq 
R_{VW}^{n}(\mathbf{1})\leq\dots \leq 
R_{VW}(\mathbf{1})\leq\mathbf{1},$$ so the limit 
$$h:=\lim_{n\rightarrow\infty}R_{VW}^n(\mathbf{1})$$
exists, and satisfies (\ref{eq6_8_2}) by the same argument.
\end{proof}
\begin{remark}\label{rem6_9}
Stated in Ruelle's thermodynamical formalism \cite{Rue89}, the 
data $\rho_W$ (measure) and $h_W$ (eigenfunction) in part (ii) of 
corollary \ref{cor6_8} represent an equilibrium distribution 
where $W$ is related to a potential function. Under mild 
conditions on $(X,r)$ and $W$, it is known that solutions 
$(\rho_W,h_W)$ exist, and we say that $\lambda_W=1$ is the 
Perron-Frobenius eigenvalue of the Ruelle operator $R_W$. In that 
case
$$1=\sup\{|\lambda|\,|\,\lambda\in\mbox{spectrum}(R_W)\}.$$
The reader is referred to \cite{Rue89}, \cite{Bal00}, 
\cite{NuLu99} and \cite{BrJo02} for further details regarding the 
spectral theory of $R_W$.
\par
Notice further that the conclusion of Ruelle's Perron-Frobenius 
theorem is a generalization of the classical Perron-Frobenius 
theorem for matrices with non-negative entries.
\end{remark}
\section{\label{extre}Extreme points}
\par
Theorem \ref{th6_no_no} shows that the map 
$\hat\mu\mapsto\mu_0:=\hat\mu\circ\theta_0^{-1}$ establishes a 
bijective affine correspondence between the following two sets:
$$M_{qi,1}^V(\xir):=\{\hat\mu\in 
M(\xir)\,|\,\hat\mu(\xir)=1,\hat\mu\mbox{ is 
}V\mbox{-quasi-invariant}\}$$ and
$$M_1^V(X):=\{\mu_0\in M(X)\,|\,\mu_0(X)=1, T_V(\mu_0)=\mu_0\}.$$
It is easy two see that both these sets are convex, and an 
application of the Banach-Alaoglu theorem \cite{Rud91} shows 
that they are compact in the weak$^*$-topology. Then, using the 
Krein-Milman theorem \cite{Rud91}, we conclude that each of these 
sets is the convex weak$^*$-closure of their extreme points. 
Moreover, since the correspondence is affine it preserves the 
extreme points. 
\par
This section is devoted to an analysis of the extreme points. 
Before we state and prove our main result on extreme points we 
need to define the concepts of conditional expectation 
$E_{\mu_0}$, and relative ergodicity.
\begin{proposition}\label{prop7_1}
Let $(X,\mathfrak{B}), r$ and $V$ be as above. Let $\mu_0$ be a 
measure in $M_1^V(X)$. Then for each bounded measurable function 
$g$ on $X$ there exists a bounded 
$r^{-1}(\mathfrak{B})$-measurable function $E_{\mu_0}(Vg)$ such 
that 
\begin{equation}\label{eq7_1_1}
\int_XV\,g\,f\circ r\,d\mu_0=\int_X E_{\mu_0}(Vg)f\circ 
r\,d\mu_0, \end{equation} for all bounded 
$\mathfrak{B}$-measurable functions $f$ on $X$. Moreover, this is 
unique up to $\mu_0\circ r^{-1}$-measure zero. 
\end{proposition}
\begin{proof}
The positive linear functional 
$$\Lambda_g\colon f\circ r\mapsto\int_XV\,g\,f\circ r\,d\mu_0$$
defines a measure on $(X,r^{-1}(\mathfrak{B}))$ which is 
absolutely continuous with respect to $\mu_0$. Indeed, if 
$\mu(E)=0$ then
$$0=\int_X\chi_E\,d\mu_0=\int_XV\,\chi_E\circ r\,d\mu_0$$
so $$\int_XV\,g\,\chi_E\circ r\,d\mu_0.$$ Therefore, by the 
Radon-Nikodym theorem, there exists some 
$r^{-1}(\mathfrak{B})$-measurable function $E_{\mu_0}(Vf)$ such 
that (\ref{eq7_1_1}) holds. Since
$$\left|\int_XE_{\mu_0}(Vg)\,f\circ 
r\,d\mu_0\right|\leq\|Vg\|_\infty\int_X|f\circ r|\,d\mu_0$$ it 
follows that $|E_{\mu_0}(Vg)|\leq\|Vg\|_\infty$ holds $\mu_0\circ 
r^{-1}$ a.e. 
\par
The uniqueness is also clear from the definition (\ref{eq7_1_1}).
\end{proof}
\begin{definition}\label{def7_1_no}
A measure $\mu_0\in M_1^V(X)$ is called {\it relatively ergodic}
with respect to $(r,V)$ if the only non-negative, bounded 
$\mathfrak{B}$-measurable functions $f$ on $X$ satisfying 
$$E_{\mu_0}(Vf)=E_{\mu_0}(V)f\circ r,\,\mbox{ pointwise }\mu_0\circ 
r^{-1}\mbox{-a.e.},$$ are the functions which are constant 
$\mu_0$-a.e.
\end{definition}

\begin{theorem}\label{th7_2}
Let $V\colon X\rightarrow[0,\infty)$ be bounded and measurable. Let\\ 
$\hat\mu\in M_{qi,1}^V(\xir)$, and 
$\mu_0:=\hat\mu\circ\theta_0^{-1}\in M_1^V(X)$. The following 
affirmations are equivalent:
\begin{enumerate}
\item $\hat\mu$ is an extreme point of $M_{qi,1}^V(\xir)$; \item 
$V\circ\theta_0\,d\hat\mu$ is ergodic with respect to $\hat r$; 
\item $\mu_0$ is an extreme point of $M_1^V(X)$; \item $\mu_0$ is 
relatively ergodic with respect to $(r,V)$.
\end{enumerate}
\end{theorem}
\begin{proof}
The arguments in the beginning of this section (mainly theorem 
\ref{th6_no_no}) show that (i) and (iii) are equivalent.
\par
We now prove (i)$\Rightarrow$(ii). Suppose $\hat\mu$ is not 
ergodic. Then there exists two measurable subset $A$ and $B$ of 
$\xir$ such that 
\begin{itemize}\item
$\int_AV\circ\theta_0\,d\hat\mu>0$ and 
$\int_BV\circ\theta_0\,d\hat\mu>0$ (note that then $\hat\mu(A)>0$ 
and $\hat\mu(B)>0$), \item $A\cup B=\xir$, and \item $A$ and $B$ 
are $\hat r$-invariant.\end{itemize} Define then the measures
$$\hat\mu_A(E)=\frac{\hat\mu(A\cap E)}{\hat\mu(A)},\quad 
\hat\mu_B(E)=\frac{\hat\mu(B\cap E)}{\hat\mu(B)}$$ for all 
$E\in\mathfrak{B}$. 
\par
Note that 
\begin{equation}\label{eq7_2_1}
\hat\mu=\hat\mu(A)\mu_a+\hat\mu(B)\mu_B.
\end{equation}
We prove next that $\mu_A,\mu_B$ are in $M_{qi,1}^V(\xir)$. 
Clearly they have total mass equal to $1$, so we only have to 
prove the $V$-quasi-invariance. For $f$ bounded, measurable 
function on $\xir$, we have:
\begin{align*}
\int_{\xir}f\circ\hat r^{-1} 
\,d\mu_A&=\frac{1}{\hat\mu(A)}\int_{\xir}f\circ\hat r^{-1}\chi_A\,d\hat\mu\\
&=\frac{1}{\hat\mu(A)}\int_{\xir}f\circ\hat r^{-1}\chi_A\circ\hat 
r^{-1} \,d\hat\mu\qquad(\hat r(A)=A)\\
&=\frac{1}{\hat\mu(A)}\int_{\xir}V\circ\theta_0\,f\,\chi_A\,d\hat\mu\\
&=\int_{\xir}V\circ\theta_0\,f\,d\mu_A.
\end{align*}
Hence $\mu_A$ is $V$-quasi-invariant. The same argument works for 
$\mu_B$. Therefore $\mu_A$ and $\mu_B$ are both in 
$M_{qi,1}^V(\xir)$. Now, equation (\ref{eq7_2_1}) contradicts the 
fact that $\hat\mu$ is  extreme. 
\par
(ii)$\Rightarrow$(i): Suppose there are some measures 
$\mu_1,\mu_2\in M_{qi,1}^V(\xir)$ and some $\lambda\in [0,1]$ 
such that 
\begin{equation}\label{eq7_2_2}
\hat\mu=\lambda\mu_1+(1-\lambda)\mu_2.
\end{equation}
Then $\mu_1$ and $\mu_2$ are absolutely continuous with respect 
to $\hat\mu$. Therefore there exist $f_1$ and $f_2$ in 
$L^1(\xir,\hat\mu)$ such that 
$$d\mu_1=f_1\,d\hat\mu,\quad d\mu_2=f_2\,d\hat\mu.$$
 Since $\mu_1$ and 
$\hat\mu$ are $V$-quasi-invariant, we have, for all bounded 
measurable functions $f$ on $\xir$,
\begin{align*}
\int_{\xir}f\,f_1\circ\hat 
r\,V\circ\theta_0\,d\hat\mu&=\int_{\xir}f\circ\hat 
r^{-1}\,f_1\,d\hat\mu\\
&=\int_{\xir}f\circ\hat r^{-1}\,d\mu_1\\
&=\int_{\xir}V\circ\theta_0\,f\,d\mu_1\\
&=\int_{\xir}V\circ\theta_0\,f\,f_1\,d\hat\mu.
\end{align*}
Therefore $f_1=f_1\circ r$ pointwise 
$V\circ\theta_0\,d\hat\mu$-almost everywhere. The hypothesis 
implies then that $f_1$ is constant 
$V\circ\theta_0\,d\hat\mu$-a.e. The same argument shows that 
$f_2$ is constant $V\circ\theta_0\,d\hat\mu$-a.e.
\par
Since $\hat\mu$ and $\mu_1$ are $V$-quasi-invariant we have also 
that 
$$1=\int_{\xir}\mathbf{1}\circ\hat r^{-1}\,d\mu_1=\int_{\xir}V\,d\mu_1=\int_{\xir}V\,f_1\,d\hat\mu.$$ 
It follows that $f_1=1$, pointwise 
$V\circ\theta_0\,d\hat\mu$-a.e. Same is true for $f_2$. Then 
$$\int_{\xir}f\,d\mu_1=\int_{\xir}V\circ\theta_0\,f\circ\hat 
r\,d\mu_1=\int_{\xir}V\circ\theta_0\,f\circ\hat 
r\,f_1\,d\hat\mu$$$$=\int_{\xir}V\, f\circ\hat 
r\,d\hat\mu=\int_{\xir}f\,d\hat\mu.$$ Hence $\mu_1=\hat\mu=\mu_2$ 
and $\hat\mu$ is extreme.
\par
(iii)$\Rightarrow$(iv): Suppose $\mu_0$ is not relatively 
ergodic. Then there exists a bounded measurable function 
$f_1\geq0$ on $X$ such that $E_{\mu_0}(Vf_1)=E_{\mu_0}(V)f_1\circ 
r$, pointwise $\mu_0\circ r^{-1}$-a.e., and $f_1$ is not constant 
$\mu_0$-a.e. We may assume that $\int_{X}f_1\,d\mu_0=1$. Define 
the measure $d\mu_1:=f_1\,d\mu_0$. We check that $\mu_1$ is in 
$M_1^V(X)$. 
\begin{gather*}
\int_Xf\,d\mu_1=\int_Xf\,f_1\,d\mu_0=\int_XV\,f\circ r f_1\circ 
r\,d\mu_0=\int_X E_{\mu_0}(V)f_1\circ r\,f\circ 
r\,d\mu_0\\=\int_X E_{\mu_0}(Vf_1)\,f\circ 
r\,d\mu_0=\int_XV\,f_1\,f\circ r\,d\mu_0=\int_XV\,f\circ 
r\,d\mu_1.
\end{gather*}
Now choose some $0<\lambda<1$ such that $\lambda 
f_1\leq1$ and define 
$$f_2:=\frac{1-\lambda f_1}{1-\lambda},\quad d\mu_2=f_2\,d\mu_0.$$
Then $E_{\mu_0}(Vf_2)=E_{\mu_0}(V)f_2\circ r$, and the same 
calculation as before shows that $T_V(\mu_2)=\mu_2$. Note also 
that $\mu_2(X)=1$. Since $\mu_0=\lambda\mu_1+(1-\lambda)\mu_2$ 
and $f_1$ is not constant, we have that $\mu_1\neq\mu_0$. It 
follows that $\mu_0$ is not an extreme point, thus contradicting 
the hypothesis and proving (iv).
\par
(iv)$\Rightarrow$(iii): Suppose 
\begin{equation}\label{eq7_2_3}
\mu_0=\lambda\mu_1+(1-\lambda)\mu_2 \end{equation}
 for some 
$\mu_1,\mu_2\in M_1^V(X)$, $\lambda\in[0,1]$. Then $\mu_1$ and 
$\mu_2$ are absolutely continuous with respect to $\mu_0$. Let 
$f_1,f_2$ the corresponding Radon-Nikodym derivatives. The 
relation (\ref{eq7_2_3}) implies that
$$1=\lambda f_1+(1-\lambda)f_2,\qquad\mu_0\mbox{-a.e.}$$
In particular $f_1$ and $f_2$ are bounded.
\par
We know that $\mu$ and $\mu_1$ have the fixed-point property. 
Then for all bounded measurable functions $f$ on $X$,
\begin{gather*}
\int_XE_{\mu_0}(V)\,f_1\circ r\,f\circ r\,d\mu_0=\int_XV\,f\circ 
r\,f_1\circ r\,d\mu_0=\int_Xf\,f_1\,d\mu_0=\int_Xf\,d\mu_1\\
=\int_XV\,f\circ r\,d\mu_1=\int_XV\,f\circ 
r\,f_1\,d\mu_0=\int_XE_{\mu_0}(Vf_1)\,f\circ r\,d\mu_0.
\end{gather*}
Therefore $E_{\mu_0}(V)f_1\circ r=E_{\mu_0}(Vf_1)$, $\mu_0\circ 
r^{-1}$-a.e.
\par
The hypothesis implies that $f_1$ is constant $\mu_0$-a.e. Since 
$\mu_1(X)=\mu_0(X)=1$, it follows that $f_1=1$ and $\mu_1=\mu_0$, 
and therefore $\mu_0$ is extreme.
\end{proof}

\begin{acknowledgements}
 This work was supported at 
the University of Iowa by a grant from the National Science 
Foundation (NSF-USA) under a Focused Research Program, 
DMS-0139473 (FRG).
\end{acknowledgements}

\end{document}